\documentclass[letterpaper,11pt]{article}

\usepackage[margin=1in]{geometry}
\usepackage{amsmath,amssymb,amsfonts,amsthm,mathtools}
\usepackage{graphicx,xcolor,xspace,enumitem,booktabs}
\usepackage{microtype,needspace}

\usepackage[size=tiny]{todonotes}
\usepackage[bookmarks,colorlinks=true,plainpages=false,
            citecolor=red,linkcolor=blue,anchorcolor=red,urlcolor=blue]{hyperref}
\hypersetup{pdftitle={Radius Allocation in Spatial Matching},
  pdfsubject={Scale-dependent reversal and optimal activation}}

\usepackage[normalem]{ulem}

\usepackage{prettyref}
\newrefformat{eq}{(\ref{#1})}
\newrefformat{sec}{Section~\ref{#1}}
\newrefformat{fig}{Figure~\ref{#1}}
\newrefformat{tab}{Table~\ref{#1}}
\newrefformat{rmk}{Remark~\ref{#1}}
\newrefformat{clm}{Claim~\ref{#1}}
\newrefformat{claim}{Claim~\ref{#1}}
\newrefformat{def}{Definition~\ref{#1}}
\newrefformat{assump}{Assumption~\ref{#1}}
\newrefformat{cor}{Corollary~\ref{#1}}
\newrefformat{lem}{Lemma~\ref{#1}}
\newrefformat{lmm}{Lemma~\ref{#1}}
\newrefformat{prop}{Proposition~\ref{#1}}
\newrefformat{app}{Appendix~\ref{#1}}
\newrefformat{thm}{Theorem~\ref{#1}}

\theoremstyle{plain}
\newtheorem{theorem}{Theorem}
\newtheorem{lemma}{Lemma}

\theoremstyle{definition}

\newtheorem{remark}{Remark}

\newcommand{\reals}{\mathbb{R}}

\newcommand{\hypercube}[1]{H_{#1}}
\newcommand{\diff}{\mathrm{d}}
\newcommand{\ones}{\mathbf{1}}

\newcommand{\indc}[1]{\mathbf{1}_{\{#1\}}}

\newcommand{\Prob}{\mathbb{P}}
\newcommand{\Expect}{\mathbb{E}}

\newcommand{\Var}{\mathsf{Var}}

\newcommand{\norm}[1]{\left\lVert#1\right\rVert}

\DeclareMathOperator*{\argmax}{arg\,max}

\DeclareMathOperator{\OPT}{OPT}
\DeclareMathOperator{\Vol}{Vol}

\newcommand{\dimension}{k}
\newcommand{\driver}{s}
\newcommand{\rider}{d}
\newcommand{\service}{r}
\newcommand{\drivervec}{\mathsf{S}}
\newcommand{\ridervec}{\mathsf{D}}
\newcommand{\servicevec}{\mathsf{R}}
\newcommand{\comparisonvec}{\mathsf{Q}}
\newcommand{\sfX}{\mathsf{X}}

\newcommand{\A}{\mathsf{A}}

\title{Distance flexibility in spatial matching: the value of concentration}
\author{
    Taha Ameen\thanks{
        T. Ameen is with the Computer Science Department, Carnegie Mellon University, Pittsburgh PA, USA, \texttt{tameen@cmu.edu}. 
    }~
    and 
    Sophie H.\ Yu\thanks{
        S.\ H.\ Yu is with the Operations, Information and Decisions Department,   the Wharton School of Business, University of Pennsylvania, Philadelphia PA, USA,  \texttt{hysophie@wharton.upenn.edu}.
    }
}
\date{}

\begin{document}
\maketitle

\begin{abstract}
    In spatial matching markets, a supply unit’s flexibility is measured by its service radius, the maximum distance at which it can serve demand. In dimensions $k \geq 2$, we study how a platform should allocate service radii among the supply nodes subject to a budget on their sum. The platform makes this choice before observing supply and demand locations, with the objective of maximizing the expected fulfilled demand. We show that the shape of a preferred allocation depends on the total budget: under suitable conditions, large budgets favor allocations that are more uniform in the sense of majorization, while small budgets favor concentration. We also characterize a non-uniform allocation that is asymptotically optimal for a very-sparse regime, and show that the uniform allocation is suboptimal in this regime. Our results provide theoretical explanations for the radius allocation questions raised by the numerical experiments in~\cite{ameen2026uniformity}.
\end{abstract}


\begingroup
\small
\tableofcontents
\endgroup

\section{Introduction}\label{sec:introduction}

In spatial matching markets, a supply node can serve only demand that is within reach. For example, in ride-hailing platforms such as Uber, Lyft, or DiDi, a driver must be sufficiently close to a rider to provide service. For platforms with autonomous vehicles such as Waymo or Zoox, a vehicle with more energy available for pickup travel can reach riders farther away, corresponding to a larger feasible service radius. Such supply nodes are compatible with more riders, but longer pickup trips consume more energy and can increase pickup times~\cite{yang2020optimizing,wang2022ondemand}. Similar constraints arise in drone delivery, where payload, energy use, and battery capacity limit the delivery range of a drone~\cite{dorling2017vehicle,stolaroff2018energy}. 
Across these applications, a supply node's maximum service distance measures its \emph{flexibility} and determines which demand nodes it can possibly serve.

Platforms can influence this flexibility before deployment, for example through decisions about the travel energy available to individual vehicles. A related allocation problem is studied in \cite{ameen2026uniformity}, where flexibility is measured by the \emph{volume} of each supply node's service region and the budget constrains total service volume. Their uniformity principle establishes conditions under which more uniform allocations of service volume yield larger expected matchings. In feature-based matching, service volume measures the breadth of demand characteristics that a supply node can accommodate. 
In applications involving physical travel, such as ride-hailing, a natural measure of flexibility is instead the maximum \emph{distance} a supply node can cover. In dimensions two and higher, service volume grows nonlinearly with radius, so fixing total service volume and fixing total service radius impose different allocation constraints. The simulations in \cite{ameen2026uniformity} illustrate this distinction: uniform allocation can be suboptimal when the budget constrains the sum of service radii.  A theoretical understanding of this radius-allocation problem remains missing. 

To illustrate how a radius budget can arise, consider the problem of allocating travel energy across a fleet of drones before their next deployment. Suppose that energy for nontravel requirements has already been reserved and that all drones consume the same amount of energy per unit of travel distance. In this setting, each drone's achievable service radius is proportional to its allocated travel energy. An aggregate energy budget therefore imposes a constraint on the sum of service radii.

This example illustrates a resource allocation decision made before matching. Extending one supply node's service radius leaves less of the shared resource available to extend the service radius of others. We study this decision through a stylized model in which a platform allocates service radii under a budget on their sum, before observing supply and demand locations. Its objective is to maximize the expected number of matches. The central question is: \emph{Should the platform distribute this radius budget uniformly across supply nodes or concentrate it on a smaller group?}

We study a market with $n$ supply nodes and $n$ demand nodes, independently and uniformly located in a unit cube $[0,1]^k$, in dimension $k\geq2$\footnote{If $k = 1$, the radius is proportional to the volume of the service region, and the results of~\cite{ameen2026uniformity} apply.}. Subject to the budget constraint, the platform allocates the service radii among the supply nodes before either set of locations is observed. After observing the locations, it forms a maximum matching among compatible pairs. Each node can participate in at most one match, and the performance metric is the expected maximum matching size. The independent uniform locations provide a benchmark for spatial uncertainty. Within this benchmark, the model isolates how geometry and unit service capacity affect the value of flexibility.

When compatible pairs are scarce, creating more of them can improve matching. Given a fixed total radius budget, concentrating the budget on fewer supply nodes can increase the expected number of compatible pairs because nominal service volume grows more than proportionally with radius in dimensions $k\geq2$. As radii increase, supply nodes become more likely to have compatible demand, but each can still complete only one match. Concentrating the budget on a small group then limits how many supply nodes can participate in matching. An allocation that gives substantial reach to more supply nodes can therefore become preferable.

Our first result shows that the ranking of radius allocations can reverse as the radius budget increases. Consider two allocations with the same total budget: one distributes it uniformly across all supply nodes, while the other distributes it equally across a fixed fraction of them. Concentration yields more expected matches at sufficiently small budgets, whereas uniform allocation yields more at sufficiently large budgets. Both comparisons hold for all sufficiently large markets.

Our second result quantifies the value of concentration by characterizing an asymptotically optimal allocation when the radius budget is scarce. Specifically, we study a regime where the total budget of radius (appropriately normalized) grows with the market size, but does so more slowly than the number of supply nodes. In this regime, we show that a simple concentration policy achieves the optimal expected matching size to first order. The policy assigns a common positive normalized radius to a vanishing fraction of supply nodes. In contrast, spreading the same budget uniformly across all supply nodes yields a vanishing fraction of the optimal expected matching size. These results provide a theoretical explanation for the simulation findings in~\cite{ameen2026uniformity} and show how the resource used to expand compatibility can change the preferred design of a matching system.

\subsection{Related literature}\label{sec:related-literature}

\paragraph{Spatial matching under geometric constraints.}
Spatial service-platform models study pricing and supply response~\cite{bimpikis2019spatial}, dynamic dispatch \cite{ozkan2020dynamic}, and pickup-time control \cite{wang2022ondemand}.~\cite{yang2020optimizing} jointly optimize matching intervals and matching radii, whereas~\cite{kanoria2025dynamic} study how the achievable matching distance depends on market size and spatial dimension when supply or demand arrives dynamically.
Routing and energy models for drone delivery also incorporate constraints on physical reach~\cite{dorling2017vehicle,stolaroff2018energy}. Finally,~\cite{ameen2026uniformity} study the ex ante allocation of spatial service regions under a budget on their total \textit{service volume}. We instead study the ex ante allocation of an aggregate radius budget across unit-capacity supply nodes. Because service volume is nonlinear in radius in dimensions $\dimension\ge2$, reallocating a fixed radius budget can change both the distribution of reach across supply nodes and the total nominal service volume.

\paragraph{Flexibility design.}
The vast literature on process flexibility asks how to distribute compatibility across resources. For example,~\cite{jordan1995principles} show how limited flexibility can recover much of the value of full flexibility, whereas~\cite{tsitsiklis2017flexible} design sparse compatibility graphs for flexible queueing systems, and~\cite{li2026optimality} establish asymptotically optimal sparsity of random regular graphs in process-flexibility and middle-mile transportation models. In service platforms, workforce composition and supply node availability also determine flexibility~\cite{dong2020managing,lobel2024frontiers}.
More closely related to our work, \cite{freund2026two} study how to allocate a flexibility budget between the two sides of a matching platform and identify regimes favoring concentration on one side or investment in both sides. \cite{ameen2026flexibility} derive exact asymptotic matching rates for this random-graph model and sharpen the comparisons between these allocations. 

\paragraph{Matching on random geometric graphs.}
Random geometric graphs model compatibility through spatial proximity~\cite{penrose2003random}.  \cite{bordenave2013matchings} characterize limiting matching proportions through local graph limits, with explicit formulas in several locally tree-like settings. However, their explicit formulas for Galton--Watson tree limits do not apply directly to our spatial graphs, whose local neighborhoods contain short cycles.
\cite{sentenac2025online} analyze offline and online matching on an interval with a common compatibility radius.

For heterogeneous service ranges, \cite{ameen2026uniformity} characterize the asymptotic expected maximum matching size in a one-dimensional model with two service ranges through an embedded Markov chain, and derive closed-form expressions in special cases. This  result quantifies, in these special cases, the matching performance for a given volume allocation. However, it does not address which allocation is optimal under a radius constraint. This will be the focus of our present work.

\subsection{Problem setup}\label{sec:model}

Fix a dimension $\dimension\geq 1$. Denote the vectors of supply and demand locations as 
\(
    \drivervec=(\driver_1,\ldots,\driver_n)
\)
and 
\( 
    \ridervec=(\rider_1,\ldots,\rider_n),
\)
respectively. All locations are mutually independent and uniform on the unit hypercube $\hypercube{\dimension} = [0,1]^{\dimension}$. Each supply node $i$ has a \emph{service radius} $r_i$ which captures its flexibility and determines its compatibility with demand nodes. Specifically, supply $i$ may only serve demand $j$ if 
\begin{equation}\label{eq:compatibility-rule}
    \norm{\driver_i-\rider_j}_2 \leq \service_i \, n^{-1/\dimension}.
\end{equation}
Thus $\service_i$ is a \emph{normalized radius}, and $\service_i n^{-1/\dimension}$ is its \emph{physical radius}. Its nominal service volume (ignoring boundary effects) is
\(
    \kappa_{\dimension} \, \service_i^{\dimension}/n,
\)
where
\(
    \kappa_{\dimension}
\)
is the volume of the Euclidean unit ball, i.e.
\[
    \kappa_\dimension:=\frac{\pi^{\dimension/2}}{\Gamma(1+\dimension/2)} \, .
\]
Here, $\Gamma(\cdot)$ denotes the Gamma function, and we use the convention $\kappa_0=1$. The normalization ensures that a supply node with fixed positive $\service_i$ has an expected number of compatible demand nodes of constant order, even as $n$ grows.\footnote{Without the factor $n^{-1/\dimension}$, a common fixed positive service radius would give each supply node an expected number of compatible demand nodes of order $n$. In that setting, a perfect matching exists with probability tending to one as $n\to\infty$, and hence the expected fraction of matched supply nodes tends to one; see~\cite[Theorem~1.3 and Proposition~5.1]{goel2005monotone}.}

Before the locations of supply and demand nodes are observed, the platform chooses a deterministic allocation of service radii,
\(
    \servicevec=(\service_1,\ldots,\service_n)\in[0,\infty)^n .
\)
This choice is subject to a constraint on the total radius $\sum_i r_i$, which we refer to as the normalized radius \emph{budget}. 
Given $\servicevec$, let $G_n(\servicevec)$ be the bipartite graph defined by~\eqref{eq:compatibility-rule}, and let $M(G)$ denote the cardinality of a maximum matching in $G$. We measure performance by the expected maximum matching size,
\begin{equation}\label{eq:objective}
    \mu_n(\servicevec):=\Expect \big[ M(G_n(\servicevec)) \big] \, .
\end{equation}
For an average normalized radius $b \geq 0$, the total normalized radius budget is $n b$, corresponding to a total physical radius of $b \, n^{1-1/\dimension}$. The optimal expected matching size under this budget is
\begin{equation}
  \OPT_n(b):=
  \max_{\substack{\servicevec\geq 0 \\
                   \sum_{i=1}^{n} \! \service_i \, = \, n b }}
  \mu_n(\servicevec).
\end{equation}
The maximum is attained by continuity and compactness
(see Appendix~\ref{app:continuity} for details). 

\subsection{Main results}\label{sec:results}

Let $\mathbf{x} = (x_1, x_2, \cdots, x_n)$ denote a vector in $\reals^n$, and denote by $\mathbf{x}^{\downarrow}$ its decreasing rearrangement. We say that $\mathbf{x}$ majorizes $\mathbf{y}$, denoted $\mathbf{x}\succeq\mathbf{y}$, if
\[
    \sum_{i=1}^{\ell}x_i^\downarrow\ge\sum_{i=1}^{\ell}y_i^\downarrow
    ~~~~ (1\leq \ell < n),
    ~~~~~~
    \sum_{i=1}^n x_i=\sum_{i=1}^n y_i.
\]
Thus $\mathbf{x}$ is the more concentrated allocation, while $\mathbf{y}$ is the more uniform allocation. For $t>0$, denote by $t \mathbf{x}$ the vector $(tx_1,\cdots, tx_n)$. The allocations $t\servicevec_n$ and $t\comparisonvec_n$ have the same budget whenever $\servicevec_n\succeq\comparisonvec_n$. Varying $t$ changes that common budget while preserving the shape of both allocations. Under a mild condition on the allocations, our first result establishes that in all dimensions larger than one, uniformity of the radius vector (in the sense of majorization) hurts at small budgets, but helps at large budgets.

\begin{theorem}[Reversal]\label{thm:reversal-new}
    Fix $\dimension \geq 2$. Let $\servicevec_n=(\service_{n,1},\ldots,\service_{n,n})$ and $\comparisonvec_n=(q_{n,1},\ldots,q_{n,n})$ belong to $[0,L]^n$, where $0<L<\infty$ is fixed. Suppose $\servicevec_n\succeq\comparisonvec_n$, and for some $a\in(0,L]$,
    \begin{align}\label{eq:support-separation}
        \limsup_{n\to\infty}
        \frac{\#\{i:\service_{n,i}>0\}}{n}
        <
        \liminf_{n\to\infty}
        \frac{\#\{i:q_{n,i}\ge a\}}{n}.
    \end{align}
    Then there exist $0<t_{\mathrm{s}}<t_{\ell}<\infty$ such that, for each fixed $t>0$,
    \[
     \begin{array}{ll}
     \mu_n(t\servicevec_n)>\mu_n(t\comparisonvec_n),
           &0<t<t_{\mathrm{s}},\\[2pt]
     \mu_n(t\servicevec_n)<\mu_n(t\comparisonvec_n),
           &t>t_{\ell},
     \end{array}
    \]
    for all sufficiently large $n$.
\end{theorem}

\begin{remark}[Thresholds and crossing] \label{rmk:crossing}
    For each $\dimension \geq 2$, explicit choices of the thresholds $t_{\mathrm{s}}$ and $t_{\ell}$ are given in~\eqref{eq:small-threshold} and~\eqref{eq:large-threshold}.
    The condition~\eqref{eq:support-separation} compares number of supply nodes that receive any positive radius under $\servicevec_n$, with supply nodes that receive at least $a$ under $\comparisonvec_n$. The fixed cutoff $a>0$ excludes radii that vanish with $n$. The two comparisons also imply a crossing: there are fixed $0<u<v<\infty$ such that every sufficiently large market has some $t_n\in(u,v)$ with $\mu_n(t_n\servicevec_n)=\mu_n(t_n\comparisonvec_n)$. The theorem makes no assertion about uniqueness of the crossing. 
\end{remark}

\begin{remark}[Effect of dimension]
    When the dimension $\dimension = 1$, the service volume of a supply node is proportional to its service radius, and so the radius-budget and volume-budget formulations are equivalent. The volume formulation has been studied by~\cite{ameen2026uniformity}: under the assumptions of Theorem~\ref{thm:reversal-new}, the uniformity principle of~\cite{ameen2026uniformity} implies that, for every $t>0$,
    \[
        \mu_n \big( t \, \servicevec_n \big) 
        <
        \mu_n \big( t \, \comparisonvec_n \big)
    \]
    for all sufficiently large $n$. Thus, the more uniform allocation remains preferable independent of the total budget, and the above reversal theorem only applies to dimensions $\dimension \geq 2$.
\end{remark}

Our second result characterizes an allocation that is asymptotically optimal when (i) the average normalized radius budget tends to zero, and (ii) the total normalized budget grows without bound. To motivate the allocation, we first consider the probability that a supply node reaches at least one demand: for a normalized radius $\service \geq 0$, define
\[
    g_{\dimension}(\service):=
    1 - \exp\big( -\kappa_{\dimension}\service^{\dimension} \, \big).
\]
In the limit as $n\to\infty$, this is the probability that a supply node with normalized radius $\service$ has at least one compatible demand. When all supply nodes have this radius, it is also the expected fraction of non-isolated supply nodes. For
$\dimension\geq 2$, define
\begin{equation}\label{eq:efficient-radius}
    C_{\dimension}:=
    \max_{\service>0}\frac{g_{\dimension}(\service)}{\service},
    ~~~~~~
    \service_{\dimension}:=
    \argmax_{\service>0}\frac{g_{\dimension}(\service)}{\service}.
\end{equation}
For $\dimension=1$, define $C_1:=\sup_{r>0}g_1(r)/r=2$. This supremum is approached as $r\downarrow0$ and is not attained at any positive radius. For $\dimension\ge2$, the maximizer $\service_{\dimension}$ is unique (see Section~\ref{sec:optimal-upper}). It gives the greatest limiting probability of reaching at least one demand per unit of radius budget, and $C_{\dimension}$ is the corresponding maximum ratio.  
The radius $\service_{\dimension}$ suggests assigning this radius to as many supply nodes as the budget permits. We compare the resulting activation allocation with the uniform allocation, formally defined below.
\begin{itemize}
    \item \textit{The activation allocation:} For $\dimension\ge2$ and $0\le b\le \service_{\dimension}$, define the activation allocation $\A_n(b)$ as follows. Let
    \[
        m=\left\lfloor\frac{nb}{\service_{\dimension}}\right\rfloor,\qquad \eta=nb-m\service_{\dimension}.
    \]
    Assign radius $\service_{\dimension}$ to $m$ supply nodes, radius $\eta$ to one further supply node if $\eta>0$, and zero to the others. If $m=n$, then $\eta=0$. The chosen indices of supply nodes with positive radius are fixed before their locations are sampled. This definition applies for all sufficiently large $n$.
    \item  \textit{The uniform allocation:} For $b\geq0$, let
    \[
        U_n(b):=\mu_n(b\ones)
    \]
    denote the expected maximum matching size under the uniform allocation, which assigns normalized radius $b$ to every supply node. Here, $\ones$ is the all-ones vector in $\reals^n$.
\end{itemize}

\begin{theorem}[Optimal allocation in the very-sparse regime]\label{thm:optimal}
    For any sequence $(b_n)_{n\geq 1}$ satisfying $ b_n \to 0$ and $n b_n \to \infty$, we have 
    \begin{equation}\label{eq:optimal-asymptotic}
        \begin{aligned}
            \lim_{n\to\infty}
            \frac{\OPT_n(b_n)}{n \, b_n}
            &= C_{\dimension} && (\dimension\ge1),\\
            \lim_{n\to\infty}
            \frac{\mu_n(\A_n(b_n))}{n \, b_n}
            &= C_{\dimension} && (\dimension\ge2),
        \end{aligned}
    \end{equation}
    and 
    \begin{align} \label{eq:U_b}
        \lim_{n\to\infty}
        \frac{U_n(b_n)}{n \, b_n^{\dimension}}
        = \kappa_{\dimension}.
    \end{align}
    Consequently,
    \begin{equation}\label{eq:uniform-ratio}
        \lim_{n\to\infty}
        \frac{U_n(b_n)}{\OPT_n(b_n)}
        = \frac{\kappa_{\dimension}} {C_{\dimension}} \lim_{n\to\infty} b_n^{\dimension-1}  =  \indc{\dimension = 1}.
    \end{equation}
\end{theorem}

\begin{remark}
    For $\dimension\geq 2$, Theorem~\ref{thm:optimal} identifies an optimal way to allocate service radii, in the very-sparse regime. The activation allocation assigns the fixed normalized radius $\service_{\dimension}$ to approximately $nb_n/\service_{\dimension}$ supply nodes. The number of activated nodes tends to infinity, but their fraction tends to zero. Activation achieves the optimal matching value to first order, whereas uniform allocation captures a vanishing fraction of that value. 

    For $\dimension=1$, note that $\kappa_{\dimension} = C_{\dimension} = 2$. In this setting, the uniform allocation is itself asymptotically optimal under the same budget conditions, i.e.
    \(
        \lim_{n\to\infty}
            U_n(b_n)/\OPT_n(b_n)=1.
    \)
    In one dimension, every feasible allocation has at most $2nb_n$ compatible pairs in expectation, which also bounds its expected matching size. The uniform allocation attains this upper bound to first order, since $\lim_{n\to\infty} U_n(b_n)/(nb_n)= \kappa_1=2$.
\end{remark}

\paragraph{Discussion and comparison to~\cite{ameen2026uniformity}.}
\label{sec:uniformity-comparison}

Our results address the radius-allocation question raised by the numerical experiments in \cite{ameen2026uniformity}. Their uniformity theorem compares allocations that have the same total service volume, which is proportional to $\sum_i\service_i^{\dimension}$. Our budget instead fixes $\sum_i\service_i$. For $\dimension\ge2$, making two unequal radii more equal reduces their total service volume. The resulting volume allocations therefore have different sums and cannot be compared directly by their uniformity theorem.

The proof in \cite{ameen2026uniformity} derives quantitative bounds for the change in matching size, when service volume is transferred between a pair of supply nodes ($T$-transforms). Repeatedly applying this comparison allows them to compare the matching size between two sufficiently separated, bounded allocations with same total service volume. Specifically, they show that in any dimension, the more uniform allocation (in the sense of majorization) yields the larger expected matching size. The key step in the proof of~\cite{ameen2026uniformity} is to establish strong concavity of an appropriate matching function, which shows that there are  diminishing marginal returns to expanding service volume. However, this idea does not transfer to our setting: Changing the variable to radius introduces a factor proportional to $\service^{\dimension-1}$ in their derivatives. For $\dimension\geq 2$, this increasing factor can offset the diminishing marginal coverage per unit of volume. Simulation results in~\cite{ameen2026uniformity} also recognize that the uniform allocation need not be optimal for the radius formulation, and the authors do not provide any theoretical analysis of the latter formulation.

Our two theorems provide complementary results for the radius allocation. Theorem~\ref{thm:reversal-new} compares allocation sequences under the condition~\eqref{eq:support-separation}, at each fixed radius scale $t$ as the market size $n$ grows. It establishes regimes where uniformity of the allocation can help or hurt the matching size. Theorem~\ref{thm:optimal} optimizes over all feasible allocations when $b_n\to0$ and $nb_n\to\infty$, giving an optimal allocation to first order, in the very-sparse regime where competition for demand is negligible. For $\dimension\geq2$, the value of concentration in this regime can be expressed as
\[
    \frac{\mu_n(\A_n(b_n))}{U_n(b_n)}
    \sim \frac{C_{\dimension}}{\kappa_{\dimension}}b_n^{1-\dimension}
    \to\infty.
\]
Thus, the activation allocation outperforms uniform allocation by a factor that grows without bound. This contrasts with the service-volume formulation of~\cite{ameen2026uniformity}, where the uniformity principle favors more uniform allocations under the conditions of that theorem. Characterizing allocations that are asymptotically optimal in other regimes remains an open question for $\dimension\geq 2$.

\section{Proof of the reversal theorem}\label{sec:reversal-proof}

We work under the assumptions of Theorem~\ref{thm:reversal-new}, with $\dimension\ge2$. The proof uses different comparisons at small and large radius scales: At small scales, we compare the expected matching sizes through the numbers of compatible pairs. At large scales, we use the difference between the numbers of supply nodes receiving positive radii. Denote
\[
    p_n:=\#\{i:\service_{n,i}>0\},
    ~~~~~~
    m_n := \#\{i:q_{n,i}\geq a\}.
\]
The separation condition \eqref{eq:support-separation} allows us to choose constants $\beta$ and $\alpha$ such that
\begin{equation}\label{eq:count-bounds}
 \limsup_{n\to\infty}\frac{p_n}{n}
 <\beta<\alpha<
 \liminf_{n\to\infty}\frac{m_n}{n}.
\end{equation}
Thus, $0<\beta<\alpha<1$, and, for all sufficiently large $n$, it holds that
\(
 1 
 \leq p_n\le\beta n<\alpha n\le m_n \,,
\)
where the first inequality follows because the two allocations have equal total radii and $\comparisonvec_n$ has total radius at least $m_na>0$. Define
\[
    A_{\dimension}:=
    a^{\dimension}
    \left(
        \frac{\alpha^{\dimension}}{\beta^{\dimension-1}}-\alpha
    \right)>0 \, ,
    ~~~~~~
    \delta:=\alpha-\beta>0 \,.
\]
Finally, for a vector $\mathbf{x}$, define the quantities 
\[
    S_\ell(\mathbf{x}):=\sum_i x_i^\ell \, , 
    ~~~~~~
    m_\ell(\mathbf{x}):=\frac{S_\ell(\mathbf{x})}{n} \,.
\]

\subsection{Concentration wins at small budgets}
\label{sec:small-proof}

We first show that the more concentrated allocation has a uniformly larger $\dimension$-th power sum. This quantity determines the leading term in the expected number of compatible pairs.

\begin{lemma}[Power-sum gap]\label{lem:power-gap}
    Suppose $\mathbf{x},\mathbf{y}\in[0,\infty)^n$ have nonincreasing coordinates and satisfy $\mathbf{x}\succeq\mathbf{y}$. If $x_i=0$ for $i>p$ and $y_i\ge a$ for $i\le m$, where $1 \leq p < m$ and $a>0$, then
    \begin{align*}
        S_{\dimension}(\mathbf{x})-S_{\dimension}(\mathbf{y})
        \ge
        a^{\dimension}
        \left(
            \frac{m^{\dimension}}{p^{\dimension-1}}-m
        \right).
    \end{align*}
\end{lemma}

Apply Lemma~\ref{lem:power-gap} to the decreasing rearrangements of $\servicevec_n$ and $\comparisonvec_n$. On $0<u<v$, the function $v^{\dimension}/u^{\dimension-1}-v$ is decreasing in $u$ and increasing in $v$. Combining this observation with~\eqref{eq:count-bounds} gives

\begin{equation}\label{eq:profile-gap}
    \frac{
        S_{\dimension}(\servicevec_n)
        -S_{\dimension}(\comparisonvec_n)
    }{n}
    \geq A_{\dimension} \, ,
\end{equation}
for all sufficiently large $n$. We next translate this power-sum gap into a difference in expected matching sizes. The following estimate bounds the losses caused by truncation of service regions at the boundary and by compatible pairs that share a node.

\begin{lemma}\label{lem:sparse}
    Let $\dimension\ge1$ and $c_{\dimension} :=  2\dimension\kappa_{\dimension-1}/(\dimension+1)$. For $\servicevec\in[0,\infty)^n$ and $t>0$ satisfying $t \max_i\service_i \leq n^{1/\dimension}$,
    \begin{align}
        0
        & \, \leq \, 
        \kappa_{\dimension}t^{\dimension}
        m_{\dimension}(\servicevec)
        -\frac{\mu_n(t\servicevec)}{n}
        \, \leq \,
         c_{\dimension}t^{\dimension+1}n^{-1/\dimension}
        m_{\dimension+1}(\servicevec)
        +
        \frac{\kappa_{\dimension}^2t^{2\dimension}}{2n^2}
        \big[
            (n-2)S_{2\dimension}(\servicevec)
            +S_{\dimension}(\servicevec)^2
        \big] \, .
        \label{eq:sparse-bound}
    \end{align}
\end{lemma}

Fix $t>0$. Since both allocations have entries bounded by $L$, Lemma~\ref{lem:sparse} applies to both for all sufficiently large $n$. Moreover, $\servicevec_n$ has at most $\beta n$ positive
entries, so
\begin{align} \label{eq:moment-bounds}
     m_{2\dimension}(\servicevec_n)\le\beta L^{2\dimension},
     \qquad
     m_{\dimension}(\servicevec_n)\le\beta L^{\dimension},
     \qquad
     m_{\dimension+1}(\servicevec_n)\le L^{\dimension+1}.
\end{align}
Using the lower bound on $\mu_n(t\servicevec_n)$ and the upper bound on $\mu_n(t\comparisonvec_n)$ from \eqref{eq:sparse-bound}, we get that for all sufficiently large $n$,
\begin{align*}
    \frac{\mu_n(t \, \servicevec_n)-\mu_n(t \, \comparisonvec_n)}{n}
    &\geq \kappa_{\dimension}t^{\dimension}
        \big[m_{\dimension}(\servicevec_n)-m_{\dimension}(\comparisonvec_n)\big]
        -c_{\dimension}t^{\dimension+1}n^{-1/\dimension}
         m_{\dimension+1}(\servicevec_n)\\
    & ~~~~~~ 
        -
        \frac{\kappa_{\dimension}^2t^{2\dimension}}{2}
        \left[
            \left(1-\frac{2}{n}\right)m_{2\dimension}(\servicevec_n)
            +m_{\dimension}(\servicevec_n)^2
        \right]\\
    &\geq \kappa_{\dimension}A_{\dimension}t^{\dimension}
        -c_{\dimension}t^{\dimension+1}n^{-1/\dimension}L^{\dimension+1}
        -\frac{\kappa_{\dimension}^2t^{2\dimension}}{2}
         \beta(1+\beta)L^{2\dimension}.
\end{align*}
Here, the last inequality substitutes~\eqref{eq:profile-gap},~\eqref{eq:moment-bounds} and uses $0\leq 1-2/n \leq 1$ for $n\geq 2$. The first term on the right-hand side of~\eqref{eq:sparse-bound} is at most $c_{\dimension}t^{\dimension+1}n^{-1/\dimension}L^{\dimension+1}$, which tends to zero because $t$ is fixed. Therefore, taking the lower limit yields
\begin{align}
    \liminf_{n\to\infty}
    \frac{
        \mu_n(t\servicevec_n)-\mu_n(t\comparisonvec_n)
    }{n}
    \geq 
    \kappa_{\dimension}A_{\dimension}t^{\dimension}
    -
    \frac{\kappa_{\dimension}^2}{2}
    \beta(1+\beta)L^{2\dimension}t^{2\dimension}.
    \label{eq:small-gap}
\end{align}
The positive term is of order $t^{\dimension}$, whereas the
remaining error is of order $t^{2\dimension}$. Define
\begin{equation}\label{eq:small-threshold}
 t_{\mathrm{s}}:=
 \left(
   \frac{2A_{\dimension}}
   {\kappa_{\dimension}\beta(1+\beta)L^{2\dimension}}
 \right)^{1/\dimension}.
\end{equation}
For every fixed $0<t<t_{\mathrm{s}}$, the right-hand side of~\eqref{eq:small-gap} is strictly positive. Hence $\mu_n(t\servicevec_n)>\mu_n(t\comparisonvec_n)$ for all sufficiently large $n$, proving the result whenever $t < t_{\mathrm{s}}$.

\subsection{Uniformity wins at large budgets}
\label{sec:large-proof}

A zero-radius supply node is isolated almost surely, regardless of the value of $t$. Therefore,
\begin{equation}\label{eq:support-upper}
    \mu_n(t\servicevec_n)\leq p_n \leq \beta n
    \qquad\text{for all sufficiently large }n.
\end{equation}
To obtain a larger matching under $\comparisonvec_n$, we match supply and demand within spatial cells. The cells are chosen small enough that every selected supply node can reach any demand node in the same cell. The next estimate bounds the matching loss caused by imbalances between the numbers of supply and demand
nodes across cells.

\begin{lemma}[Grid matching]\label{lem:grid}
    If $\servicevec$ has at least $m\in\{0,\cdots,n\}$ entries that are at least $h>0$, then
    \begin{align*}
        \mu_n(\servicevec)\ge
        m-\sqrt{\frac{m(J_n(h)^{\dimension}-1)}{2}},
        ~~~~~~ \mbox{where }
        J_n(h):=
        \left\lceil
            \frac{\sqrt{\dimension}\,n^{1/\dimension}}{h}
        \right\rceil.
    \end{align*}
\end{lemma}
Fix $t>0$ and let $\ell_n:=\lfloor\alpha n\rfloor$. For all sufficiently large $n$, at least $\ell_n$ entries of $t \, \comparisonvec_n$ are at least $ta$. Lemma~\ref{lem:grid} therefore gives
\[
     \frac{\mu_n(t\comparisonvec_n)}{n}
     \geq
     \frac{\ell_n}{n}
     -
     \sqrt{
       \frac12\frac{\ell_n}{n}
       \frac{J_n(ta)^{\dimension}-1}{n}
     } \, .
\]
Note also that
\[
  \lim_{n\to\infty}\frac{\ell_n}{n}=\alpha,
  ~~~~~~
  \lim_{n\to\infty}
  \frac{J_n(ta)^{\dimension}-1}{n}
  =
  \frac{\dimension^{\dimension/2}}{(ta)^{\dimension}}.
\]
Combining this lower bound with \eqref{eq:support-upper} yields
\begin{equation}\label{eq:large-gap}
     \liminf_{n\to\infty}
     \frac{
       \mu_n(t\comparisonvec_n)-\mu_n(t\servicevec_n)
     }{n}
     \geq
     \delta-
     \frac{
       \dimension^{\dimension/4}\sqrt{\alpha}
     }{
       \sqrt{2}\,(ta)^{\dimension/2}
     }.
\end{equation}
The error on the right decreases with $t$, while the gap $\delta$ remains fixed. Define
\begin{equation}\label{eq:large-threshold}
 t_{\ell}:=
 \max\left\{
   2t_{\mathrm{s}},
   \left(
     \frac{\dimension^{\dimension/2}\alpha}
     {2\delta^2a^{\dimension}}
   \right)^{1/\dimension}
 \right\}.
\end{equation}
For every fixed $t>t_{\ell}$, the right-hand side of~\eqref{eq:large-gap} is strictly positive. Hence $\mu_n(t\comparisonvec_n)>\mu_n(t\servicevec_n)$ for all sufficiently large $n$. This proves the result when $t>t_{\ell}$ and completes the proof of Theorem~\ref{thm:reversal-new}.

\subsection{The crossing}\label{sec:crossing-proof}

We now verify the crossing assertion stated in Remark~\ref{rmk:crossing}. Let
\[
    u:=t_{\mathrm{s}}/2,\qquad
    v:=2t_{\ell},\qquad
    f_n(t):=
        \mu_n(t\servicevec_n)-\mu_n(t\comparisonvec_n).
\]
The two comparisons imply that
\(
    f_n(u)>0>f_n(v)
\)
for all sufficiently large $n$. It remains to show that $f_n$ is continuous. For a fixed nonnegative vector $\sfX$ and scales $t,\tau\ge0$, coupling the graphs on the same locations gives
\begin{equation}\label{eq:scale-continuity}
 \left|
   \mu_n(t\sfX)-\mu_n(\tau\sfX)
 \right|
 \leq
 \kappa_{\dimension} \, 
 \big|t^{\dimension}-\tau^{\dimension}\big| \, 
 S_{\dimension}(\sfX).
\end{equation}

This bound, proved in Appendix~\ref{app:continuity}, establishes continuity of each term in $f_n$. The intermediate value theorem therefore gives a scale $t_n\in(u,v)$ satisfying $f_n(t_n)=0$ for every sufficiently large $n$.

\section{Proof of the optimal allocation theorem}\label{sec:optimal-proof}

Fix $\dimension\ge2$. We first derive an upper bound on the optimal matching value from the expected number of non-isolated supply nodes. We then show that activation attains this bound to first order. Finally, in Section~\ref{sec:uniform-proof}, we compare the uniform allocation with the optimum for every $\dimension\ge1$, treating $\dimension=1$ separately.

\subsection{An upper bound for every allocation}
\label{sec:optimal-upper}

We begin by verifying that the radius $\service_{\dimension}$ is well defined and establishing a linear bound on $g_{\dimension}(\service)$.
For $\service>0$ and $z=\kappa_{\dimension}\service^{\dimension}$, define $f(z):=e^z-1-\dimension z$. Then
\[
     \frac{\diff}{\diff\service}
    \left(\frac{g_{\dimension}(\service)}{\service}\right)
    =
    \frac{(1+\dimension z)e^{-z}-1}{\service^2}
    =-\frac{e^{-z}}{\service^2}f(z).
\]
The function $f$ satisfies
\[
     f(0)=0,\quad
    f'(0)=1-\dimension<0,\quad
    f''(z)=e^z>0,\quad
    \lim_{z\to\infty}f(z)=\infty.
\]
Strict convexity therefore implies that $f$ has exactly one positive zero, denoted $z_{\dimension}$. The derivative of $g_{\dimension}(\service)/\service$ has the opposite sign to $f(z)$, so it is positive for $0<z<z_{\dimension}$ and negative for $z>z_{\dimension}$. Thus the ratio attains its maximum at a unique positive radius. By \eqref{eq:efficient-radius}, the maximizing radius is 
\(
    \service_{\dimension}
    =(z_{\dimension}/\kappa_{\dimension})^{1/\dimension},
\)
where $z_{\dimension}$ is the unique positive solution of $e^z=1+\dimension z$, and 
\begin{align}\label{eq:scalar-bound}
  g_{\dimension}(\service)\le C_{\dimension}\service
  ~~~~~~
  \text{for all }\service\geq 0 \, .
\end{align}

Let $N_S^+(G)$ denote the number of non-isolated supply nodes in $G$. Since every matched supply node must have positive degree, $M(G)\leq N_S^+(G)$. The following lemma gives an upper bound valid for every allocation and a lower bound for a subset of supply nodes with a common radius. 

\begin{lemma}\label{lem:activation}
    For every $\servicevec\in[0,\infty)^n$,
    \begin{equation}\label{eq:nonisolation-bound}
        \mu_n(\servicevec)
        \leq
         \Expect \big[N_S^+(G_n(\servicevec)) \big]
        \leq
         \sum_{i=1}^n g_{\dimension}(\service_i)+1.
    \end{equation}
    For $m\in\{0,\ldots,n\}$ and $\service\geq0$, let $G_{n,m}(\service)$ denote the graph induced by a fixed set of $m$ supply nodes, each with normalized radius $\service$, and all $n$ demand nodes. If $\service n^{-1/\dimension}\leq 1/2$, then
    \begin{equation}\label{eq:activation-lower}
        \Expect[M(G_{n,m}(\service))]
        \geq
         m g_{\dimension}(\service)
        -2 \, \dimension\service m n^{-1/\dimension}
        -\frac{\kappa_{\dimension}^2\service^{2\dimension}}{2n}
        m(m-1).
    \end{equation}
\end{lemma}

For any budget $B\geq 0$, equations~\eqref{eq:scalar-bound} and~\eqref{eq:nonisolation-bound} imply
\begin{equation}\label{eq:opt-upper}
    \OPT_n(B/n)
    \leq \max_{\substack{\servicevec\ge0\\\sum_i\service_i=B}}
    \left(\sum_i g_{\dimension}(\service_i)+1\right)
    \le C_{\dimension}B+1.
\end{equation}
This upper bound requires no restriction on individual radii beyond the total budget constraint.

\subsection{Analysis of the activation allocation}\label{sec:activation-proof}

Denote $B_n=nb_n$ and $m_n=\lfloor B_n/\service_{\dimension}\rfloor$. Since $b_n\to0$, for all sufficiently large $n$ we have $b_n\leq \service_{\dimension}$ and $\service_{\dimension}n^{-1/\dimension}\leq 1/2$. Thus the activation allocation is well defined and the lower bound in Lemma~\ref{lem:activation} applies.
For a lower bound on its performance, restrict the compatibility graph to the $m_n$ supply nodes assigned radius $\service_{\dimension}$ and all $n$ demand nodes. This restriction discards the supply node with residual radius if one
is present, since discarding it cannot increase the maximum matching size. Applying Lemma~\ref{lem:activation} gives
\begin{align}
    \mu_n(\A_n(b_n))
    &\geq
     m_n g_{\dimension}(\service_{\dimension})
    -2 \, \dimension\service_{\dimension}m_n n^{-1/\dimension}
    -\frac{\kappa_{\dimension}^2
             \service_{\dimension}^{2\dimension}}{2n}
    m_n(m_n-1). \label{eq:LB}
\end{align}
We now express each term on the right in terms of the budget $B_n$. Since $m_n=\lfloor B_n/\service_{\dimension}\rfloor$,
\[
    \frac{B_n}{\service_{\dimension}}-1
    \le m_n\le\frac{B_n}{\service_{\dimension}}.
\]
Using the lower bound on $m_n$, together with $g_{\dimension}(\service_{\dimension})=C_{\dimension}\service_{\dimension}$ and $g_{\dimension}(\service_{\dimension})\leq 1$, gives
\begin{align*}
     m_n \, g_{\dimension}(\service_{\dimension})
    \geq
     \left(\frac{B_n}{\service_{\dimension}}-1\right)
     g_{\dimension}(\service_{\dimension})
     =
     C_{\dimension}B_n-g_{\dimension}(\service_{\dimension})
    \geq
    C_{\dimension}B_n-1.
\end{align*}
For the two remaining terms in~\eqref{eq:LB}, the upper bound on $m_n$ gives
\begin{align*}
    2\dimension\service_{\dimension}m_n n^{-1/\dimension}
    &\leq 
    2 \, \dimension B_n n^{-1/\dimension},\\
    \frac{\kappa_{\dimension}^2\service_{\dimension}^{2\dimension}}{2n}
    m_n(m_n-1)
    &
    \leq \frac{\kappa_{\dimension}^2\service_{\dimension}^{2\dimension}}{2n}m_n^2
    \le\frac{\kappa_{\dimension}^2\service_{\dimension}^{2\dimension-2}}{2}
       \frac{B_n^2}{n}.
\end{align*}
Substituting in~\eqref{eq:LB} yields
\[
    \mu_n(\A_n(b_n))
    \ge C_{\dimension}B_n-1
    -2\dimension B_n n^{-1/\dimension}
    -\frac{\kappa_{\dimension}^2\service_{\dimension}^{2\dimension-2}}{2}
     \frac{B_n^2}{n}.
\]
Since activation uses total budget $B_n$, it is feasible and satisfies $\mu_n(\A_n(b_n)) \leq \OPT_n(b_n)$. Applying~\eqref{eq:opt-upper} with $B=B_n$ also gives $\OPT_n(b_n) \leq C_{\dimension}B_n+1$.

\medskip
\noindent To prove the limits in~\eqref{eq:optimal-asymptotic}, we divide these bounds by $B_n=nb_n$, which is positive for all sufficiently large $n$. We obtain
\begin{align} \label{eq:UB-LB}
     C_{\dimension}
     -\frac{1}{B_n}
     -2\dimension n^{-1/\dimension}
     -\frac{\kappa_{\dimension}^2
                 \service_{\dimension}^{2\dimension-2}}{2}b_n
    \leq
     \frac{\mu_n(\A_n(b_n))}{B_n}
     \leq
     \frac{\OPT_n(b_n)}{B_n}
     \leq C_{\dimension}+\frac{1}{B_n}.
\end{align}
Here $\dimension$, $\kappa_{\dimension}$, and $\service_{\dimension}$ are fixed, while
\(
    \frac1{B_n}\to0,
     n^{-1/\dimension}\to 0,
\)
and
\(
    b_n\to0.
\)
The lower and upper bounds in~\eqref{eq:UB-LB} therefore both converge to $C_{\dimension}$. Since both $\mu_n(\A_n(b_n))/B_n$ and $\OPT_n(b_n)/B_n$ lie between these bounds, they also converge to $C_{\dimension}$, proving~\eqref{eq:optimal-asymptotic}.

\subsection{Analysis of the uniform allocation}\label{sec:uniform-proof}

To prove~\eqref{eq:U_b}, we bound $\kappa_{\dimension}nb_n^{\dimension}-U_n(b_n)$ and show that its ratio to $nb_n^{\dimension}$ tends to zero. Apply Lemma~\ref{lem:sparse} with $\servicevec=\ones$ and $t=b_n$. For this vector, $m_{\dimension}(\ones)=m_{\dimension+1}(\ones)=1$ and $S_{\dimension}(\ones)=S_{2\dimension}(\ones)=n$. The lemma's condition $b_n\le n^{1/\dimension}$ holds for all sufficiently large $n$, since $b_n\to0$. Substituting these values into~\eqref{eq:sparse-bound} and multiplying by $n$ gives
\begin{align*}
    0\leq
     \kappa_{\dimension}n b_n^{\dimension}-U_n(b_n)
    &\leq
     c_{\dimension}n^{1-1/\dimension}b_n^{\dimension+1}
    +\frac{\kappa_{\dimension}^2b_n^{2\dimension}}{2n}
     \big[(n-2)n+n^2\big]\\
    &=
     c_{\dimension}n^{1-1/\dimension}b_n^{\dimension+1}
    +\kappa_{\dimension}^2(n-1)b_n^{2\dimension}.
\end{align*}
Dividing by $nb_n^{\dimension}>0$ and using $(n-1)/n\le1$, we obtain
\begin{equation}\label{eq:uniform-limit}
    0\leq
     \kappa_{\dimension}
    -\frac{U_n(b_n)}{nb_n^{\dimension}}
    \leq
     c_{\dimension}n^{-1/\dimension}b_n
    +\kappa_{\dimension}^2b_n^{\dimension}
     \to 0 \, .
\end{equation}
Both terms on the right tend to zero because $b_n\to0$, proving~\eqref{eq:U_b} for every $\dimension\ge1$. We now compare uniform allocation with the optimum separately for $\dimension\ge2$ and $\dimension=1$.

\medskip 
\noindent For $\dimension\ge2$, equations~\eqref{eq:uniform-limit} and~\eqref{eq:optimal-asymptotic} give
\[
    \frac{U_n(b_n)}
      {\OPT_n(b_n)}
    =
    b_n^{\dimension-1}
    \frac{U_n(b_n)/(nb_n^{\dimension})}
      {\OPT_n(b_n)/(nb_n)}
    \to 0,
\]
since the quotient on the right converges to $\kappa_{\dimension}/C_{\dimension}$ and $b_n^{\dimension-1}\to0$.

\medskip
\noindent For $\dimension=1$, the service interval of supply $i$ has length at most $2\service_i/n$, which also bounds its compatibility probability with any fixed demand node. Therefore every feasible allocation satisfies
\[
    \mu_n(\servicevec)
    \leq \Expect \big[|E(G_n(\servicevec))|\big]
    \leq \sum_{i=1}^n n\frac{2\service_i}{n}
    =2nb_n.
\]
Since the uniform allocation is feasible, it follows that
\(
    {U_n(b_n)}/{(nb_n)}
    \leq
    {\OPT_n(b_n)}/{(nb_n)}
    \leq
    2.
\)
By~\eqref{eq:uniform-limit}, the left-hand side tends to $\kappa_1=2$. Thus $\OPT_n(b_n)/(nb_n)\to2=C_1$, proving the optimal-value limit in~\eqref{eq:optimal-asymptotic} for $\dimension=1$, and
\[
    \frac{U_n(b_n)}{\OPT_n(b_n)}
    =\frac{U_n(b_n)/(nb_n)}{\OPT_n(b_n)/(nb_n)}
    \to\frac22=1.
\]
This establishes~\eqref{eq:uniform-ratio} in both cases.
This completes the proof of Theorem~\ref{thm:optimal}.

\bibliographystyle{alpha}
\bibliography{bibliography}

\clearpage
\appendix

\section{Supplemental proofs for the reversal theorem}\label{app:reversal}

\subsection{Continuity}\label{app:continuity}

Fix $\mathbf{x},\mathbf{y}\in[0,\infty)^n$ and use the same sampled supply and demand locations to construct $G_n(\mathbf{x})$ and $G_n(\mathbf{y})$. Denote their edge sets by $E_{\mathbf{x}}$ and $E_{\mathbf{y}}$ respectively. Further, let $A\triangle B:=(A\setminus B)\cup(B\setminus A)$ denote the symmetric difference of the sets $A$ and $B$, so $E_{\mathbf{x}}\triangle E_{\mathbf{y}}$ consists of the edges present in exactly one graph.

Deleting the edges in $E_{\mathbf{x}}\setminus E_{\mathbf{y}}$ from a maximum matching in $G_n(\mathbf{x})$ leaves a matching in $G_n(\mathbf{y})$. Hence $M(G_n(\mathbf{x}))-M(G_n(\mathbf{y}))\le |E_{\mathbf{x}}\setminus E_{\mathbf{y}}|$. Interchanging $\mathbf{x}$ and $\mathbf{y}$ gives the reverse bound, and so
\begin{align*}
    |\mu_n(\mathbf{x})-\mu_n(\mathbf{y})|
    \leq \Expect\big[|M(G_n(\mathbf{x}))-M(G_n(\mathbf{y}))|\big]
    \leq \Expect\big[|E_{\mathbf{x}}\triangle E_{\mathbf{y}}|\big].
\end{align*}
For a fixed supply--demand pair $(i,j)$, the edge changes only if $\rider_j$ lies in the annulus between physical radii $x_i n^{-1/\dimension}$ and $y_i n^{-1/\dimension}$ centered at $\driver_i$. Conditional on $\driver_i$, the demand location is uniform on the unit cube, so this probability is bounded by the volume of the full annulus, namely $\kappa_{\dimension}|x_i^{\dimension}-y_i^{\dimension}|/n$. This bound does not depend on $\driver_i$, so it also holds without conditioning. Summing over all supply--demand pairs gives
\begin{align*}
    \Expect\big[|E_{\mathbf{x}}\triangle
     E_{\mathbf{y}}|\big]
    = \sum_{i=1}^n\sum_{j=1}^n
    \Prob\big((i,j)\in E_{\mathbf{x}}\triangle E_{\mathbf{y}}\big)
    \leq \sum_{i=1}^n\sum_{j=1}^n
    \frac{\kappa_{\dimension}}{n}|x_i^{\dimension}-y_i^{\dimension}| =\kappa_{\dimension}\sum_{i=1}^n|x_i^{\dimension}-y_i^{\dimension}|.
\end{align*}
Combining these inequalities proves continuity of $\mu_n$ on $[0,\infty)^n$, since the final bound tends to zero as $\mathbf{y}\to\mathbf{x}$. For $\sfX\in[0,\infty)^n$ and $t,\tau\ge0$, taking $\mathbf{x}=t\sfX$ and $\mathbf{y}=\tau\sfX$, and using $\sum_i|t^{\dimension}X_i^{\dimension}-\tau^{\dimension}X_i^{\dimension}|=|t^{\dimension}-\tau^{\dimension}|S_{\dimension}(\sfX)$, gives~\eqref{eq:scale-continuity}.
Finally, for every $b\ge0$, the feasible set
\[
    \Big\{\servicevec\in[0,\infty)^n:
        \sum_{i=1}^n\service_i=nb\Big\}
\]
is nonempty (it contains $b\ones$), closed, and bounded, since each coordinate lies in $[0,nb]$. It is therefore compact, and the continuous function $\mu_n$ attains its maximum on this set. Thus the maximum defining $\OPT_n(b)$ exists.

\subsection{Proof of Lemma~\ref{lem:power-gap}}\label{app:power-gap}

Let $\Delta_i=x_i-y_i$. For $u,v\ge0$, define
\[
    D_{\dimension}(u,v):=u^{\dimension}-v^{\dimension}-\dimension v^{\dimension-1}(u-v).
\]
Convexity of $z\mapsto z^{\dimension}$ gives $u^{\dimension}\ge v^{\dimension}+\dimension v^{\dimension-1}(u-v)$, so $D_{\dimension}(u,v)\ge0$. By definition,
\[
     S_{\dimension}(\mathbf{x})-S_{\dimension}(\mathbf{y})
     =\dimension\sum_{i=1}^n y_i^{\dimension-1}\Delta_i
      +\sum_{i=1}^n D_{\dimension}(x_i,y_i).
\]
To bound the power-sum difference from below by $\sum_{i=1}^n D_{\dimension}(x_i,y_i)$, we first show that $\sum_{i=1}^n y_i^{\dimension-1}\Delta_i\ge0$. Set $T_j=\sum_{i=1}^j\Delta_i$, with $T_0=0$. Majorization gives $T_j\ge0$ for $1\le j<n$ and $T_n=0$. Substituting $\Delta_i=T_i-T_{i-1}$ and collecting the coefficient of each $T_j$, we obtain
\begin{align*}
     \sum_{i=1}^n y_i^{\dimension-1}\Delta_i
     & =\sum_{i=1}^n y_i^{\dimension-1}(T_i-T_{i-1})\\
     & =y_n^{\dimension-1}T_n-y_1^{\dimension-1}T_0
       +\sum_{j=1}^{n-1}T_j(y_j^{\dimension-1}-y_{j+1}^{\dimension-1})\\
     & =\sum_{j=1}^{n-1}T_j(y_j^{\dimension-1}-y_{j+1}^{\dimension-1})\ge0.
\end{align*}
The last equality uses $T_0=T_n=0$, and the inequality follows from $T_j\ge0$ and $y_j\ge y_{j+1}$. Substituting into the decomposition above gives
\[
     S_{\dimension}(\mathbf{x})-S_{\dimension}(\mathbf{y})
     \geq \sum_{i=1}^n D_{\dimension}(x_i,y_i).
\]
For $p<i \leq m$, we have $x_i=0$ and $y_i\geq a$, so
\[
     D_{\dimension}(x_i,y_i)
     =(\dimension-1)y_i^{\dimension}
     \geq (\dimension-1)a^{\dimension} \, .
\]
Keeping the terms with $i\le m$ and using $D_{\dimension}(x_i,y_i)\ge0$ for $i>m$, we obtain
\begin{equation}\label{eq:remainder-gap}
    \begin{aligned}
        S_{\dimension}(\mathbf{x})-S_{\dimension}(\mathbf{y})
        &\geq \sum_{i=1}^p D_{\dimension}(x_i,y_i)
            +(\dimension-1)\sum_{i=p+1}^m y_i^{\dimension}\\
        &\geq \sum_{i=1}^p D_{\dimension}(x_i,y_i)
            +(m-p)(\dimension-1)a^{\dimension}.
 \end{aligned}
\end{equation}

It remains to bound $\sum_{i=1}^p D_{\dimension}(x_i,y_i)$. Since $\Delta_i=-y_i\le0$ for $i>p$ and $\sum_{i=1}^n\Delta_i=0$, we have
\[
     P:= \sum_{i=1}^p(\Delta_i)_+
    = \sum_{i=1}^n(-\Delta_i)_+
    \geq \sum_{i=p+1}^m y_i\ge(m-p)a.
\]
Thus the values $(\Delta_i)_+$ among the first $p$ coordinates have sum at least $(m-p)a$. To use this bound, for $v,w\ge0$ define
\[
    F_v(w) := D_{\dimension}(v+w,v)
    =\sum_{j=2}^{\dimension}\binom{\dimension}{j}v^{\dimension-j}w^j.
\]
Note that $F_v(w)$ is nondecreasing in $v$, and that upon setting $v=a$, we have that $F_a(w)$ is convex and nondecreasing in $w$. For $i\leq p$, we have $y_i\geq a$ because $p<m$. Consequently,
\begin{align*}
    D_{\dimension}(x_i,y_i)
     &=F_{y_i}(\Delta_i)\ge F_a(\Delta_i),
     &&\text{if }\Delta_i>0,\\
     D_{\dimension}(x_i,y_i)
     &\ge0=F_a(0),
     &&\text{if }\Delta_i\le0.
\end{align*}
In both cases, $D_{\dimension}(x_i,y_i)\ge F_a((\Delta_i)_+)$. Summing over $i\leq p$ and applying Jensen's inequality gives
\begin{align*}
     \sum_{i=1}^p D_{\dimension}(x_i,y_i)
     \geq \sum_{i=1}^p F_a((\Delta_i)_+)
     \geq p \, F_a\left(\frac1p\sum_{i=1}^p(\Delta_i)_+\right)
      = p \, F_a(P/p)
     \geq p \, F_a \big((m-p)a/p \big),
\end{align*}
where the last inequality uses $P\ge(m-p)a$ and monotonicity of $F_a$.
Substituting this bound into~\eqref{eq:remainder-gap}, and using $F_a(w)=(a+w)^{\dimension}-a^{\dimension}-\dimension a^{\dimension-1}w$, gives
\begin{align*}
     S_{\dimension}(\mathbf{x})-S_{\dimension}(\mathbf{y})
     &
     \geq p \, F_a \big((m-p)a/p \big)+(m-p)(\dimension-1)a^{\dimension}\\
     &=p\left[
           \left(\frac{ma}{p}\right)^{\dimension}-a^{\dimension}-\dimension \frac{m-p}{p}a^{\dimension}\right]+(m-p)(\dimension -1)a^{\dimension}\\
     &=a^{\dimension}\left(\frac{m^{\dimension}}{p^{\dimension -1}}-m\right),
\end{align*}
which is the desired bound.

\subsection{Proof of Lemma~\ref{lem:sparse}}\label{app:sparse}

For a finite simple graph $G$, let $C(G)=\sum_{v\in V(G)}\binom{\deg_G(v)}2$. This counts the unordered pairs of distinct edges that share a vertex. For every such pair of edges in the original graph $G$, select one of its two edges, and then delete the set of all selected edges. An edge may be selected for more than one pair, so at most $C(G)$ distinct edges are deleted. No two surviving edges can share a vertex: their pair would have caused at least one of them to be selected for deletion. The surviving edges therefore form a matching. Consequently,
\begin{equation}\label{eq:edge-conflicts}
    |E(G)|-C(G) \leq M(G) \leq |E(G)| \, .
\end{equation}

Let $U,V$ be independent and uniform on $[0,1]^{\dimension}$. Their difference has density $\prod_{\ell=1}^{\dimension}(1-|z_\ell|)$ on $[-1,1]^{\dimension}$. Thus, for $0\leq h\leq 1$,
\begin{align}
    0  \leq \kappa_{\dimension}h^{\dimension}-\Prob(\norm{U-V}_2\leq h)
    &= \int_{B(0,h)}
         \left[1-\prod_{\ell=1}^{\dimension}(1-|z_\ell|)\right]\,\diff z\notag\\
    & \leq \sum_{\ell=1}^{\dimension}\int_{B(0,h)}|z_\ell|\,\diff z\notag\\
    &= 2\dimension \kappa_{\dimension -1}\int_0^h z(h^2-z^2)^{(\dimension -1)/2}\,\diff     z\notag\\
    &= c_{\dimension}h^{\dimension +1}.
    \label{eq:boundary-bound}
\end{align}
The inequality uses $1-\prod_\ell(1-u_\ell)\leq \sum_\ell u_\ell$ for $u_\ell\in[0,1]$. Apply~\eqref{eq:boundary-bound} with $h_i=t\service_i n^{-1/\dimension}$ and sum over the $n$ demands for each supply:
\begin{equation}\label{eq:edge-expectation}
     \kappa_{\dimension}t^{\dimension}S_{\dimension}(\servicevec)
     -c_{\dimension}t^{\dimension +1}n^{-1/\dimension}S_{\dimension +1}(\servicevec)
     \leq \Expect|E(G_n(t\servicevec))|
     \leq \kappa_{\dimension}t^{\dimension}S_{\dimension}(\servicevec).
\end{equation}

Set $a_i=\kappa_{\dimension}t^{\dimension}\service_i^{\dimension}/n$. Conditional on the location of supply $i$, its incident edges from distinct demands are independent, each with probability at most $a_i$. Thus the expected number of edge pairs sharing a supply is at most
\[
    \binom n2\sum_i a_i^2
    = \frac{\kappa_{\dimension}^2t^{2\dimension}(n-1)}{2n}S_{2\dimension}(\servicevec).
\]
Conditional on a demand location, its incident edges from distinct supplies
are independent, with respective probabilities at most $a_i$.
The expected number of edge pairs sharing a demand is at most
\[
    n\sum_{i<j}a_i a_j
    =\frac{\kappa_{\dimension}^2t^{2\dimension}}{2n}
   \left[S_{\dimension}(\servicevec)^2-S_{2\dimension}(\servicevec)\right].
\]
Adding gives
\begin{equation}\label{eq:conflict-bound}
    \Expect \big[ C(G_n(t\servicevec)) \big]
        \le\frac{\kappa_{\dimension}^2t^{2\dimension}}{2n}
       \left[(n-2)S_{2\dimension}(\servicevec)+S_{\dimension}(\servicevec)^2\right].
\end{equation}
Take expectations in \eqref{eq:edge-conflicts}, substitute~\eqref{eq:edge-expectation} and \eqref{eq:conflict-bound}, and divide by $n$. This proves \eqref{eq:sparse-bound}.

\subsection{Proof of Lemma~\ref{lem:grid}}\label{app:grid}

Choose $m$ supply nodes whose radii are at least $h$, and choose any $m$ demands, with both index sets fixed before sampling locations. Partition $[0,1]^{\dimension}$ into $J^{\dimension}$ equal cubes, where $J=J_n(h)$. Every cell has diameter $\sqrt{\dimension}/J\le h n^{-1/\dimension}$, so every selected supply--demand pair in a common cell is compatible.

Let $X_C,Y_C$ be the selected supply and demand counts in cell $C$. Matching within each cell gives
\[
    M(G_n(\servicevec))\ge\sum_C\min\{X_C,Y_C\}
    =m-\frac12\sum_C|X_C-Y_C|.
\]
The equality uses $\min\{x,y\}=(x+y-|x-y|)/2$ and $\sum_C X_C=\sum_C Y_C=m$. For each cell, $X_C$ and $Y_C$ are independent $\operatorname{Bin}(m,J^{-\dimension})$ random variables. Hence
\[
     \Expect|X_C-Y_C|
     \leq \sqrt{\Expect[(X_C-Y_C)^2]}
     =\sqrt{2mJ^{-\dimension}(1-J^{-\dimension})}.
\]
The inequality is Cauchy--Schwarz. For the equality, $\Expect(X_C-Y_C)=0$, and independence gives $\Expect[(X_C-Y_C)^2]=\Var(X_C)+\Var(Y_C)$; each binomial variance is $mJ^{-\dimension}(1-J^{-\dimension})$.
Summing over the $J^{\dimension}$ cells gives
\[
     \mu_n(\servicevec)
     \geq m-\frac{J^{\dimension}}{2}\sqrt{2mJ^{-\dimension}(1-J^{-\dimension})}
     =m-\sqrt{\frac{m(J^{\dimension}-1)}2}.
\]

\section{Supplemental proof for the optimal allocation theorem}\label{app:optimal}

\subsection{Proof of Lemma~\ref{lem:activation}}\label{app:activation}

\paragraph{Proof of~\eqref{eq:nonisolation-bound}.}
First, for every integer $n\ge1$ and $z\ge0$, note that
\begin{equation}\label{eq:binomial-envelope}
    0 \leq e^{-z}-(1-z/n)_+^n \leq \frac1n.
\end{equation}
The lower bound follows from $1-u\leq e^{-u}$. For $0\leq u\leq 1/2$,
\[
    \frac{\diff }{\diff u}\big(\log(1-u)+u+u^2\big)
    =\frac{u(1-2u)}{1-u}\ge0.
\]
Since $\log(1-u)+u+u^2$ vanishes at $u=0$, this gives $\log(1-u)\ge-u-u^2$ on $[0,1/2]$. Substituting $u=z/n$ and exponentiating yields $(1-z/n)^n\ge e^{-z-z^2/n}$ when $z\le n/2$.
If $z\le n/2$, it follows that
\[
     e^{-z}-(1-z/n)^n
     \leq e^{-z}(1-e^{-z^2/n})
     \leq \frac{z^2e^{-z}}n
     \leq \frac4{e^2n}
     < \frac1n.
\]
If $z>n/2$, the difference is at most $e^{-n/2}\le1/n$, since $\sup_{x>0}xe^{-x/2}=2/e<1$. This proves~\eqref{eq:binomial-envelope}.

Let $v_i=\Vol(B(\driver_i,\service_i n^{-1/\dimension})\cap [0,1]^{\dimension})$, so
$v_i\le\min\{\kappa_{\dimension}\service_i^{\dimension}/n,1\}$.
Conditioning on $\driver_i$ gives
 \begin{align*}
     \Prob(\deg(i)>0\mid \driver_i)
     &=1-(1-v_i)^n\\
     &\leq 1-(1-\kappa_{\dimension}\service_i^{\dimension}/n)_+^n
     \leq g_{\dimension}(\service_i)+\frac1n.
\end{align*}
Summing over $i$ and using that $M(G)\leq N_S^+(G)$ yields~\eqref{eq:nonisolation-bound}.

\paragraph{Proof of~\eqref{eq:activation-lower}.}
For any finite bipartite graph, select one incident edge at each nonisolated
supply. If $c_j$ selected edges meet demand $j$, keep one of them.
The resulting matching has size
\[
     \sum_j\indc{c_j>0}
     =N_S^+(G)-\sum_j(c_j-1)_+
     \geq N_S^+(G)-\sum_j\binom{\deg_G(j)}2.
\]
The equality uses $\sum_j c_j=N_S^+(G)$, since exactly one edge was selected at each nonisolated supply. The inequality follows from $(c_j-1)_+\le\binom{c_j}{2}\leq\binom{\deg_G(j)}2$ for every demand $j$. Thus
\begin{equation}\label{eq:demand-conflicts}
    M(G) \geq N_S^+(G)-\sum_j\binom{\deg_G(j)}2.
\end{equation}
For $G_{n,m}(\service)$, conditioning on a demand location leaves the $m$ supply
locations independent. Each selected supply reaches that demand with
probability at most $\kappa_{\dimension}\service^{\dimension}/n$. Therefore
\begin{equation}\label{eq:activation-conflicts}
     \Expect\left[\sum_{j=1}^n\binom{\deg(j)}2\right]
     \leq n \binom{m}{2} 
        \left( \frac{\kappa_{\dimension}\service^{\dimension}}{n}\right)^2
     =\frac{\kappa_{\dimension}^2\service^{2\dimension}}{2n}m(m-1).
\end{equation}
Let $h=\service n^{-1/\dimension}\le1/2$. A selected supply lies in $[h,1-h]^{\dimension}$ with probability $(1-2h)^{\dimension}$. On that event its radius-$h$ ball is contained in $[0,1]^{\dimension}$, so
\begin{align*}
        \Prob(\deg(i)>0)
        & \geq (1-2h)^{\dimension}\left[1-(1-\kappa_{\dimension}\service^{\dimension}/n)^n\right]\\
        & \geq (1-2h)^{\dimension}g_{\dimension}(\service)
        \geq g_{\dimension}(\service)-2\dimension h.
\end{align*}
The last inequality uses $(1-2h)^{\dimension}\geq 1-2\dimension h$ and $g_{\dimension}(\service)\leq 1$. Summing over the $m$ supplies, and applying~\eqref{eq:demand-conflicts}--\eqref{eq:activation-conflicts} gives \eqref{eq:activation-lower}. This completes the proof.

\end{document}